\documentclass[12pt,reqno,twoside]{article}

\usepackage{amsfonts}
\usepackage{amsmath}
\usepackage{amssymb}
\usepackage{amsthm}
\usepackage{enumitem}
\usepackage[top=36.9mm,bottom=20mm,left=20mm,right=20mm]{geometry}
\usepackage[pdfstartview=FitV,bookmarks=false]{hyperref}
\usepackage{mathpazo}

\def\func#1{\mathrm{#1}}

\newcommand{\lt}{\mathcal{L}}
\newcommand{\reg}[1]{\mathcal{R}^{#1}}

\newcommand{\cnt}{\mathrm{C}}
\newcommand{\creg}[1]{\mathcal{R}_{c}^{#1}}
\newcommand{\crd}[1]{\mathrm{C}_{\mathrm{rd}}^{#1}}
\newcommand{\C}{\mathbb{C}}
\newcommand{\N}{\mathbb{N}}
\newcommand{\R}{\mathbb{R}}
\newcommand{\T}{\mathbb{T}}
\newcommand{\Z}{\mathbb{Z}}

\newcommand{\constl}{L}
\newcommand{\constm}{M}
\newcommand{\constn}{N}

\newcommand{\kernel}{K}
\newcommand{\forcing}{f}

\newtheorem{theorem}{Theorem}
\newtheorem{definition}{Definition}
\newtheorem{lemma}{Lemma}
\newtheorem{example}{Example}
\theoremstyle{definition}

\newtheorem{remark}{Remark}

\numberwithin{equation}{section}

\begin{document}
\markboth{\centerline{B.~Karpuz}}{\centerline{Volterra integral equations}}
\title{Basics of Volterra integral equations\\ on time scales}
\author{Ba\c{s}ak KARPUZ\footnote{On leave at the University of Calgary, Canada.\newline\textbf{Address}: Department of Mathematics, Faculty of Science and Arts, ANS Campus, Afyon Kocatepe University, 03200 Afyonkarahisar, Turkey. \textbf{Email}: bkarpuz@gmail.com \textbf{Web}. \url{http://www2.aku.edu.tr/\string~bkarpuz}}}

\maketitle

\begin{abstract}
This paper studies existence and uniqueness of solutions to generalized Volterra integral equations.
Since our proof for existence and uniqueness does not make use of Banach fixed point theorem unlike the previous papers focused on this subject,
we can replace continuity property of the kernel function with the weaker one rd-continuity.
The paper also covers results concerning the following concepts:
The notion of resolvent kernel,
and its role in formulation of the solution,
the reciprocity property of kernels,
Piccard iterates,
relation between linear dynamic equations and Volterra integral equations,
some special type of kernels together with several illustrative examples.
\end{abstract}

\section{Introduction}\label{intro}

The so-called \emph{Volterra integral equations} find application in demography, the study of viscoelastic materials, and in insurance mathematics through the renewal equation.
Let $a\in\R$ be a fixed point and $\kernel$ be a complex-valued continuous function defined on the domain $\{(t,s)\in\R^{2}:\ b\geq t\geq s\geq a\}$.
Let $\forcing$ be a complex-valued continuous function defined on $[a,b]$, then as is well known, Volterra integral equation of the second kind and of the first kind are defined by
\begin{equation}
\varphi(t)=\lambda\int_{a}^{t}\kernel(t,\eta)\varphi(\eta)\mathrm{d}\eta+\forcing(t)\quad\text{for}\ t\in[a,b],\notag
\end{equation}
where $\lambda$ is a complex parameter, and
\begin{equation}
\int_{a}^{t}\kernel(t,\eta)\varphi(\eta)\mathrm{d}\eta=\forcing(t)\quad\text{for}\ t\in[a,b],\notag
\end{equation}
respectively.
In the above equations, $\varphi$ defined on $[a,b]$ is a complex-valued continuous function to be solved.
The readers are referred the books \cite{MR2155102,MR1105769} for the fundamentals of this huge theory.
And for the previous results on the dynamic generalization of these equations,
the readers may refer to \cite{MR2669151,MR2548121} and the references cited therein.

The paper is arranged as follows:
In \S~\ref{secpots}, we shall give a short account on the notion of time scale concept;
in \S~\ref{secgviesk}, we will introduce and start studying generalized Volterra integral equations of the second kind;
in \S~\ref{secgviefk}, generalized Volterra integral equations of the first kind will be introduced together with their relations with the second kind.
Finally, in \S~\ref{fc}, we conclude the paper by making some comments on nonlinear Volterra integral equations and systems.

\section{Preliminaries on Time Scales}\label{secpots}

For a reader not familiar with the time scale calculus, we find helpful to introduce the following introductory information.
A \emph{time scale}, which inherits the standard topology on $\R$, is a nonempty closed subset of reals.
Here, and later throughout this paper, a time scale will be denoted by the symbol $\T$, and the intervals with a subscript $\T$ are used to denote the intersection of the usual interval with $\T$.
For $t\in\T$, we define the \emph{forward jump operator} $\sigma:\T\to\T$ by $\sigma(t):=\inf(t,\infty)_{\T}$ while the \emph{backward jump operator} $\rho:\T\to\T$ is defined by $\rho(t):=\sup(-\infty,t)_{\T}$, and the \emph{graininess function} $\mu:\T\to\R_{0}^{+}$ is defined to be $\mu(t):=\sigma(t)-t$.
A point $t\in\T$ is called \emph{right-dense} if $\sigma(t)=t$ and/or equivalently $\mu(t)=0$ holds; otherwise, it is called \emph{right-scattered}, and similarly \emph{left-dense} and \emph{left-scattered} points are defined with respect to the backward jump operator.
A function $f:\T\to\R$ is said to be \emph{Hilger differentiable} (or $\Delta$-differentiable) at the point $t\in\T$ if there exists $\ell\in\R$ such that for any $\varepsilon>0$ there exists a neighborhood $U$ of $t$ such that
\begin{equation}
\big|[f(\sigma(t))-f(s)]-\ell[\sigma(t)-s]\big|\leq\varepsilon|\sigma(t)-s|\quad\text{for all}\ s\in U,\notag
\end{equation}
and this present case, we denote $f^{\Delta}(t)=\ell$.
We shall mean the Hilger derivative of a function when we only say derivative unless otherwise is specified.
A function $f$ is called \emph{rd-continuous} provided that it is continuous at right-dense points in $\T$, and has finite limit at left-dense points, and the \emph{set of rd-continuous functions} are denoted by $\crd{}(\T,\R)$.
The set of functions $\crd{1}(\T,\R)$ includes the functions $f$ whose derivative is in $\crd{}(\T,\R)$ too.
For $s,t\in\T$ and a function $f\in\crd{}(\T,\R)$, the $\Delta$-integral is defined to be
\begin{equation}
\int_{s}^{t}f(\eta)\Delta\eta=F(t)-F(s),\notag
\end{equation}
where $F\in\crd{1}(\T,\R)$ is an anti-derivative of $f$, i.e., $F^{\Delta}=f$ on $\T^{\kappa}$, where $\T^{\kappa}:=\T\backslash\{\sup\T\}$ if $\sup\T=\max\T$ and satisfies $\rho(\max\T)\neq\max\T$; otherwise, $\T^{\kappa}:=\T$.
It should be noted that $\Delta$-integral by means of the Reimann-sum is introduced in \cite{MR1942437}.

A function $f\in\crd{}(\T,\R)$ is called \emph{regressive} if $1+f\mu\neq0$ on $\T^{\kappa}$, and $f\in\crd{}(\T,\R)$ is called \emph{positively regressive} if $1+f\mu>0$ on $\T^{\kappa}$.
The \emph{set of regressive functions} and the \emph{set of positively regressive functions} are denoted by $\reg{}(\T,\R)$ and $\reg{+}(\T,\R)$, respectively, and $\reg{-}(\T,\R)$ is defined similarly.
For simplicity, we denote by $\creg{}(\T,\R)$ the \emph{set of regressive constants}, and similarly, we define the sets $\creg{+}(\T,\R)$ and $\creg{-}(\T,\R)$.

Let $f\in\reg{}(\T,\R)$, then the \emph{generalized exponential function} $\func{e}_{f}(\cdot,s)$ on a time scale $\T$ is defined to be the unique solution of the initial value problem
\begin{equation}
\begin{cases}
x^{\Delta}(t)=f(t)x(t)\quad\text{for}\ t\in\T^{\kappa}\\
x(s)=1
\end{cases}\notag
\end{equation}
for some fixed $s\in\T$.
For $h\in\R^{+}$, set $\C_{h}:=\{z\in\C:\ z\neq-1/h\}$, $\Z_{h}:=\{z\in\C:\ -\pi/h<\func{Img}(z)\leq\pi/h\}$, and $\C_{0}:=\Z_{0}:=\C$.
For $h\in\R_{0}^{+}$, define the \emph{cylinder transformation} $\xi_{h}:\C_{h}\to\Z_{h}$ by
\begin{equation}
\xi_{h}(z):=
\begin{cases}
z,&h=0\\
\displaystyle\frac{1}{h}\func{Log}(1+zh),&h>0
\end{cases}\notag
\end{equation}
for $z\in\C_{h}$, then the exponential function can also be written in the form
\begin{equation}
\func{e}_{f}(t,s):=\exp\bigg\{\int_{s}^{t}\xi_{\mu(\eta)}\big(f(\eta)\big)\Delta\eta\bigg\}\quad\text{for}\ s,t\in\T.\notag
\end{equation}
It is known that the exponential function $\func{e}_{f}(\cdot,s)$ is strictly positive on $[s,\infty)_{\T}$ provided that $f\in\reg{+}([s,\infty)_{\T},\R)$, while $\func{e}_{f}(\cdot,s)$ alternates in sign at right-scattered points of the interval $[s,\infty)_{\T}$ provided that $f\in\reg{-}([s,\infty)_{\T},\R)$.
For $h\in\R_{0}^{+}$, let $z,w\in\C_{h}$, the \emph{circle plus} and the \emph{circle minus} are respectively defined by
\begin{equation}
z\oplus_{h}w:=z+w+zwh\quad\text{and}\quad z\ominus_{h}w:=\frac{z-w}{1+wh}.\notag
\end{equation}
It is also known that $\reg{+}(\T,\R)$ is a subgroup of $\reg{}(\T,\R)$, i.e., $0\in\reg{+}(\T,\R)$, $f,g\in\reg{+}(\T,\R)$ implies $f\oplus_{\mu}g\in\reg{+}(\T,\R)$ and $\ominus_{\mu}f\in\reg{+}(\T,\R)$, where $\ominus_{\mu}f:=0\ominus_{\mu}f$ on $\T$ (see \cite{MR1843232}).

The definition of the generalized monomials on time scales $\func{h}_{k}:\T\times\T\to\R$ are given as follows:
\begin{equation}
\func{h}_{k}(t,s):=
\begin{cases}
1,&k=0\\
\displaystyle\int_{s}^{t}\func{h}_{k-1}(\eta,s)\Delta\eta,&k\in\N
\end{cases}\notag
\end{equation}
for $s,t\in\T$ (see \cite{MR1843232}).
Using induction, it is easy to see that $\func{h}_{k}(t,s)\geq0$ for all $k\in\N_{0}$ and all $s,t\in\T$ with $t\geq s$ and $(-1)^{k}\func{h}_{k}(t,s)\geq 0$ for all $k\in\N$ and all $s,t\in\T$ with $t\leq s$ (see \cite{MR2506155}).

A relation between the exponential function and the monomials is given by
\begin{equation}
\func{e}_{\lambda}(t,s)=\sum_{\ell=0}^{\infty}\lambda^{\ell}\func{h}_{\ell}(t,s)\quad\text{for}\ s,t\in\T\ \text{with}\ t\geq s,\notag
\end{equation}
where $\lambda\in\creg{}(\T,\R)$ (see \cite{MR2320804}).

For $f\in\crd{}(\T\times\T,\R)$ ($f(\cdot,s)\in\crd{}(\T,\R)$ for each fixed $s\in\T$ and $f(t,\cdot)\in\crd{}(\T,\R)$ for each fixed $t\in\T$), we have
\begin{equation}
\int_{s}^{t}\int_{s}^{\eta}f(\eta,\zeta)\Delta\zeta\Delta\eta=\int_{s}^{t}\int_{\sigma(\zeta)}^{t}f(\eta,\zeta)\Delta\eta\Delta\zeta\notag
\end{equation}
for $s\in\T^{\kappa}$ and $t\in\T$ (see \cite{MR2506155}), and at the end of the paper, we shall present a different proof of this result, which was given in \cite{MR2506155} incompletely.

Using the change of order formula above, we can obtain the alternative definition of the monomials as follows
\begin{equation}
\func{h}_{k}(t,s)=\displaystyle\int_{s}^{t}\func{h}_{k-1}\big(t,\sigma(\eta)\big)\Delta\eta\quad\text{for}\ k\in\N\notag
\end{equation}
for $s,t\in\T$ (see \cite{MR2320804}).
Therefore, for all $k\in\N$, we see that $\func{h}_{k}^{\Delta_{1}}(t,s)=\func{h}_{k-1}(t,s)$ for all $(t,s)\in\T^{\kappa}\times\T$ and $\func{h}_{k}^{\Delta_{2}}(t,s)=-\func{h}_{k-1}\big(t,\sigma(s)\big)$ for all $(t,s)\in\T\times\T^{\kappa}$.

It should be noted that if $\sum_{\ell=0}^{\infty}f_{\ell}$ is a series of $\Delta$-integrable functions on $[a,b)_{\T}$ which converges uniformly on $[a,b)_{\T}$, then
\begin{equation}
\int_{a}^{b}\sum_{\ell=0}^{\infty}f_{\ell}(\eta)\Delta\eta=\sum_{\ell=0}^{\infty}\int_{a}^{b}f_{\ell}(\eta)\Delta\eta.\notag
\end{equation}
In other words, the series can be term-by-term integrable (see \cite{MR2000143}).

The readers are referred to \cite{MR1843232} for further interesting details in the time scale theory.

\section{The Generalized Volterra Integral Equation of the Second Kind}\label{secgviesk}

Let $\T$ be a time scale and $a\in\T$ and $b\in\T$ with $b>a$.
We introduce the generalized Volterra integral equation of the second kind by
\begin{equation}
\varphi(t)=\lambda\int_{a}^{t}\kernel(t,\eta)\varphi(\eta)\Delta\eta+\forcing(t)\quad\text{for}\ t\in[a,b]_{\T}\ \text{and}\ \lambda\in\R\label{gvieskeq1}
\end{equation}
under the following primary assumptions:
\begin{enumerate}[label={(H\arabic*)},leftmargin=*,ref=(H\arabic*)]
\setcounter{enumi}{0}
\item\label{h1} $\kernel\in\crd{}(\Omega(a,b),\R)$ ($\kernel(\cdot,s)\in\crd{}([s,b]_{\T},\R)$ for each fixed $s\in[a,b]_{\T}$ and $\kernel(t,\cdot)\in\crd{}([a,t]_{\T},\R)$ for each fixed $t\in[a,b]_{\T}$), where
    \begin{equation}
    \Omega(a,b):=\{(t,s)\in\T\times\T:\ b\geq t\geq s\geq a\}.\label{gvieskeq2}
    \end{equation}
\item\label{h2} $f\in\crd{}([a,b]_{\T},\R)$.
\end{enumerate}
If we take a look at \eqref{gvieskeq1}, we see the integral element $\eta$ travels through the interval $[a,t)_{\T}\subset[a,b]_{\T}$, which means $\Omega(b,a)$ could be defined to be $\{(t,s)\in\T\times\T:\ b\geq t>s\geq a\}$.
However, this set need not be compact, i.e., $\kernel$ may not be bounded.
To over come this case, we defined $\Omega(b,a)$ as in \eqref{gvieskeq2}.

We continue with the following simple example.

\begin{example}\label{gvieskex1}
Let $n\in\N$, and consider
\begin{equation}
\varphi(t)=-\int_{a}^{t}\frac{(t-c)^{n-1}}{\sum_{k=1}^{n}\big(\sigma(\eta)-c\big)^{k}(\eta-c)^{n-k}}\varphi(\eta)\Delta\eta+\frac{(t-c)^{n-1}}{(a-c)^{n}}\quad\text{for}\ t\in[a,b]_{\T},\label{gvieskex1eq1}
\end{equation}
where $c\in\R$ with $c\not\in[a,b]_{\T}$.
For this equation, we have $\lambda=-1$, $\kernel(t,s)=(t-c)^{n-1}/\big(\sum_{k=1}^{n}\big(\sigma(\eta)-c\big)^{k}(\eta-c)^{n-k}\big)$ and $\forcing(t)=(t-c)^{n-1}/(a-c)^{n}$ for $t\in[a,b]_{\T}$.
Using \cite[Theorem~1.24(ii)]{MR1843232} and letting $\varphi(t)=1/(t-c)$ for $t\in[a,b]_{\T}$, and substituting $\varphi$ into \eqref{gvieskeq1}, we have
\begin{align}
\frac{1}{t-c}+\int_{a}^{t}\frac{(t-c)^{n-1}}{\sum_{k=1}^{n}\big(\sigma(\eta)-c\big)^{k}(\eta-c)^{n+1-k}}\Delta\eta=&\frac{1}{t-c}-\frac{(t-c)^{n-1}}{(\eta-c)^{n}}\bigg|_{\eta=a}^{\eta=t}\notag\\
=&\frac{(t-c)^{n-1}}{(a-c)^{n}}\notag
\end{align}
for all $t\in[a,b]_{\T}$.
This shows that $\varphi$ solves \eqref{gvieskex1eq1}.
\end{example}

\subsection{Existence and Uniqueness of Solutions}\label{seceauos}

In this section, we prove existence and uniqueness of solutions to \eqref{gvieskeq1}, which is the most important result of this paper.

\begin{theorem}[Existence and Uniqueness]\label{eauosthm1}
Assume that \ref{h1} and \ref{h2} hold.
Then, for each $\lambda\in\R$, Equation~\eqref{gvieskeq1} admits a unique solution $\varphi\in\crd{}([a,b]_{\T},\R)$.
\end{theorem}

\begin{proof}
Suppose for now that the formal power series of the form
\begin{equation}
\sum_{\ell=0}^{\infty}\varphi_{\ell}(t)\lambda^{\ell}\quad\text{for}\ t\in[a,b]_{\T},\label{eauosthmprfeq1}
\end{equation}
where $\{\varphi_{n}\}_{n\in\N_{0}}$ is a sequence of continuous functions defined on $[a,b]_{\T}$, satisfies \eqref{eauosthmprfeq1}.
We shall show later that the function series \eqref{eauosthmprfeq1} converges absolutely and uniformly on $[a,b]_{\T}$ for $|\lambda|<\infty$.
Substituting \eqref{eauosthmprfeq1} into \eqref{gvieskeq1}, and arranging the resultant series, we get
\begin{align}
\forcing(t)=&\sum_{\ell=0}^{\infty}\varphi_{\ell}(t)\lambda^{\ell}-\lambda\int_{a}^{t}\kernel(t,\eta)\sum_{\ell=0}^{\infty}\varphi_{\ell}(\eta)\lambda^{\ell}\Delta\eta\notag\\
=&\varphi_{0}(t)+\sum_{\ell=1}^{\infty}\bigg(\varphi_{\ell}(t)-\int_{a}^{t}\kernel(t,\eta)\varphi_{\ell-1}(\eta)\Delta\eta\bigg)\lambda^{\ell}\notag
\end{align}
for all $t\in[a,b]_{\T}$.
By a comparison of the coefficients, we get the following relations
\begin{equation}
\varphi_{n}(t)=
\begin{cases}
\forcing(t),&n=0\\
\displaystyle\int_{a}^{t}\kernel(t,\eta)\varphi_{n-1}(\eta)\Delta\eta,&n\in\N
\end{cases}\quad\text{for}\ t\in[a,b]_{\T}.\label{eauosthmprfeq2}
\end{equation}
Since $\kernel$ and $\forcing$ are rd-continuous functions, there exist two constants $\constl,\constm\in\R^{+}$ such that $|\kernel(t,s)|\leq\constm$ for all $(t,s)\in\Omega(a,b)$ and $|\forcing(t)|\leq\constl$ for all $t\in[a,b]_{\T}$.
Recall that every rd-continuous function (more truly regulated) is bounded on compact domains. 
Then, we have
\begin{equation}
|\varphi_{0}(t)|\leq\constl\quad\text{and}\quad|\varphi_{1}(t)|\leq\constl\constm\func{h}_{1}(t,a)\quad\text{for all}\ t\in[a,b]_{\T},\label{eauosthmprfeq3}
\end{equation}
and hence
\begin{align}
|\varphi_{2}(t)|\leq&\constl\constm\int_{a}^{t}|\kernel(t,\eta)|\func{h}_{1}(\eta,a)\Delta\eta\leq\constl\constm^{2}\int_{a}^{t}\func{h}_{1}(\eta,a)\Delta\eta\notag\\
=&\constl\constm^{2}\func{h}_{2}(t,a)\notag
\end{align}
for all $t\in[a,b]_{\T}$.
By the emerging pattern, we feel that
\begin{equation}
|\varphi_{n}(t)|\leq\constl\constm^{n}\func{h}_{n}(t,a)\quad\text{for all}\ t\in[a,b]_{\T}\ \text{and}\ n\in\N.\label{eauosthmprfeq4}
\end{equation}
Indeed, assuming that \eqref{eauosthmprfeq4} is true for $n\in\N$ and using \eqref{eauosthmprfeq2}, we get
\begin{align}
|\varphi_{n+1}(t)|\leq&\constl\constm^{n}\int_{a}^{t}|\kernel(t,\eta)|\func{h}_{n}(\eta,a)\Delta\eta\leq\constl\constm^{n+1}\int_{a}^{t}\func{h}_{n}(\eta,a)\Delta\eta\notag\\
=&\constl\constm^{n+1}\func{h}_{n+1}(t,a)\notag
\end{align}
for all $t\in[a,b]_{\T}$, which proves that \eqref{eauosthmprfeq4} is true for all $n\in\N$.
It follows that
\begin{equation}
\sum_{\ell=0}^{\infty}\big|\varphi_{\ell}(t)\lambda^{\ell}\big|\leq\constl\sum_{\ell=0}^{\infty}\big(|\lambda|\constm\big)^{\ell}\func{h}_{\ell}(t,a)=\constl\func{e}_{|\lambda|\constm}(t,a)\leq\constl\func{e}_{|\lambda|\constm}(b,a)\notag
\end{equation}
for all $t\in[a,b]_{\T}$ ($|\lambda|\constm\geq0$ implies $|\lambda|\constm\in\creg{+}(\T,\R)$), which proves that the series \eqref{eauosthmprfeq1} converges absolutely and uniformly on $[a,b]_{\T}$ for $|\lambda|<\infty$.
Hence, $\varphi$ defined by
\begin{equation}
\varphi(t):=\sum_{\ell=0}^{\infty}\varphi_{\ell}(t)\lambda^{\ell}\quad\text{for}\ t\in[a,b]_{\T},\label{eauosthmprfeq5}
\end{equation}
where the sequence $\{\varphi_{n}\}_{n\in\N_{0}}$ is defined recursively by \eqref{eauosthmprfeq2}, is a solution of \eqref{gvieskeq1}.
Up to here, we have showed existence of a solution, and what follows next is the proof of the uniqueness part.
Suppose that there exist two different solutions, and denote by $\psi\in\crd{}([a,b]_{\T},\R_{0}^{+})$ the absolute value of the difference of these two solutions.
Then, we see that $\psi$ satisfies
\begin{equation}
\psi(t)=\lambda\int_{a}^{t}\kernel(t,\eta)\psi(\eta)\Delta\eta\leq|\lambda|\constm\int_{a}^{t}\psi(\eta)\Delta\eta\quad\text{for all}\ t\in[a,b]_{\T}.\label{eauosthmprfeq6}
\end{equation}
Applying the Gr\"{o}nwall inequality (see \cite[Theorem~6.4]{MR1843232}) to \eqref{eauosthmprfeq6}, we see that $\psi(t)\leq0$ for all $t\in[a,b]_{\T}$, i.e., $\psi\equiv0$ on $[a,b]_{\T}$.
This proves uniqueness of the solutions, and completes the proof.
\end{proof}

\begin{remark}\label{eauosrmk1}
We have to note that if $\kernel\in\cnt{}(\Omega(a,b),\R)$ and $f\in\cnt{}([a,b]_{\T},\R)$,
then the unique solution of \eqref{gvieskeq1} is continuous.
\end{remark}

\subsection{Resolvent Kernels}\label{secrk}

In this section, we shall obtain the unique solution of \eqref{gvieskeq1} in terms of some special functions.
For this purpose, we introduce first the definition of iterated kernels.

\begin{definition}[Iterated Kernels]\label{rkdf1}
Assume that \ref{h1} holds.
The sequence of rd-continuous functions $\{\kernel_{n}\}_{n\in\N_{0}}$ defined recursively by
\begin{equation}
\kernel_{n}(t,s):=
\begin{cases}
\kernel(t,s),&n=0\\
\displaystyle\int_{\sigma(s)}^{t}\kernel(t,\eta)\kernel_{n-1}(\eta,s)\Delta\eta,&n\in\N
\end{cases}
\quad\text{for}\ (t,s)\in\Omega(a,b)\ \text{and}\ n\in\N_{0}\label{rkdf1eq1}
\end{equation}
is called the sequence \emph{iterated kernels} of \eqref{gvieskeq1}.
\end{definition}

The following result shows an interesting property of the iterated kernels.

\begin{lemma}\label{rklm1}
Assume that \ref{h1} holds.
Then the iterated kernels $\{\kernel_{n}\}_{n\in\N_{0}}$ of \eqref{gvieskeq1} satisfy
\begin{equation}
\kernel_{n}(t,s)=\int_{\sigma(s)}^{t}\kernel_{n-1}(t,\eta)\kernel(\eta,s)\Delta\eta\quad\text{for}\ (t,s)\in\Omega(a,b)\ \text{and}\ n\in\N.\label{rklm1eq1}
\end{equation}
\end{lemma}

\begin{proof}
The claim is true for $n=1$ trivially.
Assume now that it is true for some $n\in\N$.
For all $(t,s)\in\Omega(a,b)$, we have
\begin{align}
\kernel_{n+1}(t,s)=&\int_{\sigma(s)}^{t}\kernel(t,\eta)\kernel_{n}(\eta,s)\Delta\eta\notag\\
=&\int_{\sigma(s)}^{t}\kernel(t,\eta)\bigg(\int_{\sigma(s)}^{\eta}\kernel_{n-1}(\eta,\zeta)\kernel(\zeta,s)\Delta\zeta\bigg)\Delta\eta\notag\\
=&\int_{\sigma(s)}^{t}\int_{\sigma(s)}^{\eta}\kernel(t,\eta)\kernel_{n-1}(\eta,\zeta)\kernel(\zeta,s)\Delta\zeta\Delta\eta\notag\\
=&\int_{\sigma(s)}^{t}\int_{\sigma(\zeta)}^{t}\kernel(t,\eta)\kernel_{n-1}(\eta,\zeta)\kernel(\zeta,s)\Delta\zeta\Delta\eta\notag\\
=&\int_{\sigma(s)}^{t}\bigg(\int_{\sigma(\zeta)}^{t}\kernel(t,\eta)\kernel_{n-1}(\eta,\zeta)\Delta\eta\bigg)\kernel(\zeta,s)\Delta\zeta\notag\\
=&\int_{\sigma(s)}^{t}\kernel_{n}(t,\zeta)\kernel(\zeta,s)\Delta\zeta,\notag
\end{align}
where we have applied the change of order formula (see Theorem~\ref{ssecithm11}) for $\Delta$-integrals while passing to the fourth line.
The claim is hence true for $(n+1)$ too, and this completes the last step of the mathematical induction.
We have thus proved that \eqref{rklm1eq1} is true.
\end{proof}

With the following lemma, we provide an upper bound for the iterated kernels.

\begin{lemma}\label{rklm2}
Assume that \ref{h1} holds.
Then the iterated kernels $\{\kernel_{n}\}_{n\in\N_{0}}$ of \eqref{gvieskeq1} satisfy
\begin{equation}
|\kernel_{n}(t,s)|\leq\constm^{n+1}\func{h}_{n}(t,s)\quad\text{for all}\ (t,s)\in\Omega(a,b)\ \text{and}\ n\in\N,\label{rklm2eq1}
\end{equation}
where
\begin{equation}
\constm:=\sup_{(t,s)\in\Omega(a,b)}|\kernel(t,s)|.\notag
\end{equation}
\end{lemma}

\begin{proof}
Clearly, we have $|\kernel_{0}(t,s)|=|\kernel(t,s)|\leq\constm$ for all $(t,s)\in\Omega(a,b)$, and similarly, we obtain
\begin{equation}
|\kernel_{1}(t,s)|\leq\constm^{2}\int_{s}^{t}1\Delta\eta=\constm^{2}\func{h}_{1}(t,s)\quad\text{for all}\ (t,s)\in\Omega(a,b).\notag
\end{equation}
By a simple inductive argument, one has \eqref{rklm2eq1}.
Thus the proof is completed.
\end{proof}

The following lemma shows that the terms of the series \eqref{eauosthmprfeq5},
which converges to the unique solution $\varphi$ of \eqref{gvieskeq1},
can be written in terms of the iterated kernels $\kernel_{n}$ and the forcing term $\forcing$.

\begin{lemma}\label{rklm3}
Assume that \ref{h1} and \ref{h2} hold.
Then the sequence $\{\varphi_{n}\}_{n\in\N_{0}}$ defined by \eqref{eauosthmprfeq2} satisfies
\begin{equation}
\varphi_{n}(t)=\int_{a}^{t}\kernel_{n-1}(t,\eta)\forcing(\eta)\Delta\eta\quad\text{for all}\ t\in[a,b]_{\T}\ \text{and}\ n\in\N,\notag
\end{equation}
where $\{\kernel_{n}\}_{n\in\N_{0}}$ is the sequence of iterated kernels defined by \eqref{rkdf1eq1}.
\end{lemma}

\begin{proof}
We proceed by mathematical induction.
Trivially the claim holds for $n=1$.
Assuming its validity for $n\in\N$, we get
\begin{align}
\varphi_{n+1}(t)=&\int_{a}^{t}\kernel(t,\eta)\varphi_{n}(\eta)\Delta\eta\notag\\
=&\int_{a}^{t}\kernel(t,\eta)\bigg(\int_{a}^{\eta}\kernel_{n-1}(\eta,\zeta)\forcing(\zeta)\Delta\zeta\bigg)\Delta\eta\notag\\
=&\int_{a}^{t}\int_{a}^{\eta}\kernel(t,\eta)\kernel_{n-1}(\eta,\zeta)\forcing(\zeta)\Delta\zeta\Delta\eta\notag
\end{align}
for all $t\in[a,b]_{\T}$.
Applying the change of order formula to the right-hand side of the last equality above and using Definition~\ref{rkdf1} yield that
\begin{align}
\varphi_{n+1}(t)=&\int_{a}^{t}\int_{\sigma(\zeta)}^{t}\kernel(t,\eta)\kernel_{n-1}(\eta,\zeta)\forcing(\zeta)\Delta\eta\Delta\zeta\notag\\
=&\int_{a}^{t}\bigg(\int_{\sigma(\zeta)}^{t}\kernel(t,\eta)\kernel_{n-1}(\eta,\zeta)\Delta\eta\bigg)\forcing(\zeta)\Delta\zeta\notag\\
=&\int_{a}^{t}\kernel_{n}(t,\zeta)\forcing(\zeta)\Delta\zeta,\notag
\end{align}
holds for all $t\in[a,b]_{\T}$. 
This completes the proof by justifying of the last step of the mathematical induction.
\end{proof}

The notion of the resolvent kernel given below plays an important role in the formulation of the unique solution $\varphi$ of \eqref{gvieskeq1}.

\begin{definition}[The Resolvent Kernel]\label{rkdf2}
Assume that \ref{h1} holds.
The series $\Gamma:\Omega(a,b)\times\R\to\R$ defined by
\begin{equation}
\Gamma(\lambda;t,s):=\sum_{\ell=0}^{\infty}\kernel_{\ell}(t,s)\lambda^{\ell}\quad\text{for}\ (t,s)\in\Omega(a,b)\ \text{and}\ \lambda\in\R\label{rkdf3eq1}
\end{equation}
is called the \emph{resolvent kernel} of \eqref{gvieskeq1}.
\end{definition}

As the resolvent kernel $\Gamma$ is given by a series of functions, its convergence is therefore of theoretical interest.

\begin{theorem}[Convergence of the Resolvent Kernel]\label{rkthm1}
Assume that \ref{h1} holds.
Let $\kernel$ and $\Gamma$ be the kernel and the resolvent kernel of \eqref{gvieskeq1}, respectively.
Then, the resolvent kernel $\Gamma$ converges absolutely and uniformly for $|\lambda|<\infty$.
\end{theorem}

\begin{proof}
It follows from Lemma~\ref{rklm2} that
\begin{align}
\big|\Gamma(\lambda;t,s)\big|\leq&\sum_{\ell=0}^{\infty}|\kernel_{\ell}(t,s)\lambda^{\ell}|\leq\constm\sum_{\ell=0}^{\infty}|\lambda^{\ell}|\constm^{\ell}\func{h}_{\ell}(t,s)=\constm\func{e}_{|\lambda|\constm}(t,s)\notag\\
\leq&\constm\func{e}_{|\lambda|\constm}(b,a)\notag
\end{align}
for all $(t,s)\in\Omega(a,b)$ and $\lambda\in\R$ with $|\lambda|<\infty$.
An application of Weierstrass $M$-test completes the proof.
\end{proof}

An interesting property between the kernel $\kernel$ and the resolvent kernel $\Gamma$ is presented below.

\begin{theorem}[Reciprocity of Kernels]\label{rkthm2}
Assume that \ref{h1} holds.
If $\Gamma$ is the resolvent kernel of the kernel $\kernel$, then the resolvent kernel of the kernel $\Gamma$ is the kernel $\kernel$ itself.
\end{theorem}

\begin{proof}
It follows from \eqref{rkdf3eq1} that
\begin{equation}
\Gamma(\lambda;t,s)=\lambda\sum_{\ell=0}^{\infty}\kernel_{\ell+1}(t,s)\lambda^{\ell}+\kernel_{0}(t,s)\quad\text{for all}\ (t,s)\in\Omega(a,b)\ \text{and}\ \lambda\in\R.\label{rkthm2prfeq1}
\end{equation}
Applying Theorem~\ref{rkthm1}, we learn that reversal of the order of integration and the sum is permissable when $|\lambda|<\infty$.
From \eqref{rkthm2prfeq1} and Lemma~\ref{rklm1}, for all $(t,s)\in\Omega(a,b)$, we obtain
\begin{align}
\Gamma(\lambda;t,s)=&\lambda\sum_{\ell=0}^{\infty}\bigg(\int_{\sigma(s)}^{t}\kernel(t,\eta)\kernel_{\ell}(\eta,s)\Delta\eta\bigg)\lambda^{\ell}+\kernel(t,s)\notag\\
=&\lambda\int_{\sigma(s)}^{t}\kernel(t,\eta)\sum_{\ell=0}^{\infty}\kernel_{\ell}(\eta,s)\lambda^{\ell}\Delta\eta+\kernel(t,s)\notag
\end{align}
which yields
\begin{equation}
\Gamma(\lambda;t,s)=\lambda\int_{\sigma(s)}^{t}\kernel(t,\eta)\Gamma(\lambda;\eta,s)\Delta\eta+\kernel(t,s)\quad\text{for all}\ (t,s)\in\Omega(a,b).\label{rkthm2prfeq3}
\end{equation}
On the other hand, it follows from \eqref{rkthm2prfeq1} and Lemma~\ref{rklm1} that
\begin{align}
\Gamma(\lambda;t,s)=&\lambda\sum_{\ell=0}^{\infty}\bigg(\int_{\sigma(s)}^{t}\kernel_{\ell}(t,\eta)\kernel(\eta,s)\Delta\eta\bigg)\lambda^{\ell}+\kernel(t,s)\notag\\
=&\lambda\int_{\sigma(s)}^{t}\Gamma(\lambda;t,\eta)\kernel(\eta,s)\Delta\eta+\kernel(t,s)\notag
\end{align}
for all $(t,s)\in\Omega(a,b)$, or equivalently, we have
\begin{equation}
\kernel(t,s)=-\lambda\int_{\sigma(s)}^{t}\Gamma(\lambda;t,\eta)\kernel(\eta,s)\Delta\eta+\Gamma(\lambda;t,s)\quad\text{for all}\ (t,s)\in\Omega(a,b).\label{rkthm2prfeq4}
\end{equation}
Taking into account \eqref{rkthm2prfeq3} and \eqref{rkthm2prfeq4}, we complete the proof.
\end{proof}

Now, we are in a position to prove the unique solution $\varphi$ of \eqref{gvieskeq1} in terms of the resolvent kernel $\Gamma$ and the forcing term $\forcing$.

\begin{theorem}\label{rkthm3}
Assume that \ref{h1} and \ref{h2} hold.
If $\Gamma$ is the resolvent kernel of the kernel $\kernel$, then the unique solution of \eqref{gvieskeq1} is given by
\begin{equation}
\varphi(t)=\lambda\int_{a}^{t}\Gamma(\lambda;t,\eta)\forcing(\eta)\Delta\eta+\forcing(t)\quad\text{for}\ t\in[a,b]_{\T}.\notag
\end{equation}
\end{theorem}

\begin{proof}
By \eqref{eauosthmprfeq1}, \eqref{eauosthmprfeq2}, \eqref{eauosthmprfeq5} and Lemma~\ref{rklm3}, we have
\begin{align}
\varphi(t)=&\lambda\sum_{\ell=0}^{\infty}\varphi_{\ell+1}(t)\lambda^{\ell}+\varphi_{0}(t)\notag\\
=&\lambda\sum_{\ell=0}^{\infty}\bigg(\int_{a}^{t}\kernel_{\ell}(t,\eta)\forcing(\eta)\Delta\eta\bigg)\lambda^{\ell}+\forcing(t)\notag
\end{align}
for all $t\in[a,b]_{\T}$.
By Theorem~\ref{rkthm1}, we can reverse the order of integration and the sum is permissable when $|\lambda|<\infty$.
Thus, we obtain 
\begin{align}
\varphi(t)=&\lambda\int_{a}^{t}\bigg(\sum_{\ell=0}^{\infty}\kernel_{\ell}(t,\eta)\lambda^{\ell}\bigg)\forcing(\eta)\Delta\eta+\forcing(t)\notag\\
=&\lambda\int_{a}^{t}\Gamma(\lambda;t,\eta)\forcing(\eta)\Delta\eta+\forcing(t)\notag
\end{align}
for all $t\in[a,b]_{\T}$.
This completes the proof.
\end{proof}

\subsection{Picard Iterates}\label{secsa}

We start this section with a famous sequence of functions, which converges the solution uniformly.

\begin{definition}[Picard Iterates]\label{sadf1}
Assume that \ref{h1} and \ref{h2} hold.
The sequence of rd-continuous functions $\{\varphi_{n}\}_{n\in\N_{0}}$ defined by
\begin{equation}
\varphi_{n}(t):=\lambda\displaystyle\int_{a}^{t}\kernel(t,\eta)\varphi_{n-1}(\eta)\Delta\eta+\forcing(t)\quad\text{for}\ t\in[a,b]_{\T}\ \text{and}\ n\in\N,
\label{sadf1eq2}
\end{equation}
where $\varphi_{0}\in\crd{}([a,b]_{\T},\R)$ is chosen an arbitrarily, is called the sequence of \emph{Picard iterates} to \eqref{gvieskeq1}.
\end{definition}

Now, we prove another existence and uniqueness result for \eqref{gvieskeq1} under weaker conditions.

\begin{theorem}\label{sathm1}
Assume that \ref{h1} and \ref{h2} hold, and that $\{\varphi_{n}\}_{n\in\N_{0}}$ be the sequence of Picard iterates of \eqref{gvieskeq1}
Then, as $n\to\infty$ the function $\varphi_{n}$ tends uniformly to the rd-continuous unique solution of \eqref{gvieskeq1}.
\end{theorem}

\begin{proof}
We shall first show by mathematical induction that
\begin{equation}
|\varphi_{n}(t)-\varphi_{n-1}(t)|\leq|\lambda|^{n}\constl\constm^{n}\func{h}_{n}(t,a)+|\lambda|^{n-1}\constn\constm^{n-1}\func{h}_{n-1}(t,a)\quad\text{for all}\ t\in[a,b]_{\T}\label{sathm1prfeq1}
\end{equation}
and all $n\in\N$, where
\begin{equation}
\constl:=\sup_{t\in[a,b]_{\T}}|\varphi_{0}(t)|,\quad\constm:=\sup_{(t,s)\in\Omega(a,b)}|\kernel(t,s)|\quad\text{and}\quad\constn:=\sup_{t\in[a,b]_{\T}}|\forcing(t)-\varphi_{0}(t)|.\notag
\end{equation}
Clearly, we have
\begin{align}
|\varphi_{1}(t)-\varphi_{0}(t)|=&\bigg|\lambda\int_{a}^{t}\kernel(t,\eta)\varphi_{0}(\eta)\Delta\eta+\forcing(t)-\varphi_{0}(t)\bigg|\notag\\
\leq&|\lambda|\bigg|\int_{a}^{t}\kernel(t,\eta)\varphi_{0}(\eta)\Delta\eta\bigg|+\big|\forcing(t)-\varphi_{0}(t)\big|\notag\\
\leq&|\lambda|\constl\constm\func{h}_{1}(t,a)\constn\notag
\end{align}
for all $t\in[a,b]_{\T}$, which shows the claim is valid for $n=1$.
Assume now that the claim is true for some $n\in\N$, then we have
\begin{align}
|\varphi_{n+1}(t)-\varphi_{n}(t)|=&|\lambda|\bigg|\int_{a}^{t}\kernel(t,\eta)\varphi_{n}(\eta)\Delta\eta-\int_{a}^{t}\kernel(t,\eta)\varphi_{n-1}(\eta)\Delta\eta\bigg|\notag\\
=&|\lambda|\bigg|\int_{a}^{t}\kernel(t,\eta)\big[\varphi_{n}(\eta)-\varphi_{n-1}(\eta)\big]\Delta\eta\bigg|\notag\\
\leq&|\lambda|\constm\bigg|\int_{a}^{t}\Big[|\lambda|^{n}\constl\constm^{n}\func{h}_{n}(\eta,a)+|\lambda|^{n-1}\constn\constm^{n-1}\func{h}_{n-1}(\eta,a)\Big]\Delta\eta\bigg|\notag\\
=&|\lambda|^{n+1}\constl\constm^{n+1}\func{h}_{n+1}(t,a)+|\lambda|^{n}\constn\constm^{n}\func{h}_{n}(t,a)\notag
\end{align}
for all $t\in[a,b]_{\T}$, which proves that the claim is also true when $n$ is replaced with $(n+1)$.
Hence, we have justified \eqref{sathm1prfeq1}.
We shall show now that the sequence
\begin{equation}
\bigg\{\varphi_{0}+\sum_{\ell=0}^{n-1}\big[\varphi_{\ell+1}-\varphi_{\ell}\big]\bigg\}_{n\in\N_{0}}=\{\varphi_{n}\}_{n\in\N_{0}},\label{sathm1prfeq2}
\end{equation}
where the empty sum is assumed to be zero, converges uniformly.
From \eqref{sathm1prfeq1} and \eqref{sathm1prfeq2}, we have
\begin{align}
|\varphi_{n}(t)|\leq&|\varphi_{0}(t)|+\sum_{\ell=0}^{n-1}\big|\varphi_{\ell+1}(t)-\varphi_{\ell}(t)\big|\notag\\
\leq&|\varphi_{0}(t)|+\sum_{\ell=0}^{n-1}\Big[|\lambda|^{\ell+1}\constl\constm^{\ell+1}\func{h}_{\ell+1}(t,a)+|\lambda|^{\ell}\constn\constm^{\ell}\func{h}_{\ell}(t,a)\Big]\notag\\
\end{align}
for all $t\in[a,b]_{\T}$ and $n\in\N_{0}$.
Letting $n\to\infty$, we can estimate for all $t\in[a,b]_{\T}$ that
\begin{equation}
\sum_{\ell=0}^{n-1}|\lambda|^{\ell+1}\constm^{\ell+1}\func{h}_{\ell+1}(t,a)\to\func{e}_{|\lambda|\constm}(t,a)-1\quad\text{and}\quad\sum_{\ell=0}^{n-1}|\lambda|^{\ell}\constm^{\ell}\func{h}_{\ell}(t,a)\to\func{e}_{|\lambda|\constm}(t,a),\notag
\end{equation}
which implies
\begin{equation}
\sum_{\ell=0}^{\infty}\Big[|\lambda|^{\ell+1}\constl\constm^{\ell+1}\func{h}_{\ell+1}(t,a)+|\lambda|^{\ell}\constn\constm^{\ell}\func{h}_{\ell}(t,a)\Big]=(\constl+\constn)\func{e}_{|\lambda|\constm}(t,a)-\constl\leq(\constl+\constn)\func{e}_{|\lambda|\constm}(b,a)-\constl\notag
\end{equation}
for all $t\in[a,b]_{\T}$.
It follows from the Weierstrass $M$-test that the infinite series
\begin{equation}
\varphi_{0}(t)+\sum_{\ell=0}^{\infty}\big[\varphi_{\ell+1}(t)-\varphi_{\ell}(t)\big]\quad\text{for}\ t\in[a,b]_{\T},\notag
\end{equation}
whose sequence of partial sums is $\{\varphi_{n}\}_{n\in\N_{0}}$, converges uniformly.
Let $\varphi:=\lim_{n\to\infty}\varphi_{n}$ on $[a,b]_{\T}$.
Letting $n\to\infty$ on both sides of \eqref{sadf1eq2} and reversing the order of limit and the integral, we see that the function $\varphi$ satisfies \eqref{gvieskeq1}.
This together with Theorem~\ref{eauosthm1} completes the proof.
\end{proof}

Let us illustrate this result with an example. 

\begin{example}\label{saex1}
Consider the following Volterra integral equation
\begin{equation}
\varphi(t)=\int_{a}^{t}\varphi(\eta)\Delta\eta+1\quad\text{for}\ t\in[a,b]_{\T}.\label{saex1eq1}
\end{equation}
Clearly, when we compare \eqref{gvieskeq1} with \eqref{saex1eq1}, we see that $\kernel\equiv1$ on $\Omega(a,b)$ and $\forcing\equiv1$ on $[a,b]_{\T}$.
Let us compute the sequence of Picard iterates to \eqref{saex1eq1} with the initial term $\varphi_{0}=1$ on $[a,b]_{\T}$.
Then, we have
\begin{align}
\varphi_{1}(t)=&\int_{a}^{t}\Delta\eta+1=\func{h}_{1}(t,a)+1\notag\\
\varphi_{2}(t)=&\int_{a}^{t}\big(\func{h}_{1}(\eta,a)+1\big)\Delta\eta+1=\func{h}_{2}(t,a)+\func{h}_{1}(t,a)+1\notag
\end{align}
for $t\in[a,b]_{\T}$.
By repeating in the same manner, we get
\begin{equation}
\varphi_{n}(t)=\sum_{\ell=0}^{n}\func{h}_{\ell}(t,a)\quad\text{for}\ t\in[a,b]_{\T},\notag
\end{equation}
which yields
\begin{equation}
\lim_{n\to\infty}\varphi_{n}(t)=\lim_{n\to\infty}\sum_{\ell=0}^{n}\func{h}_{\ell}(t,a)=\func{e}_{1}(t,a)\quad\text{for}\ t\in[a,b]_{\T}.\notag
\end{equation}
This is exactly the unique solution of \eqref{saex1eq1}.
\end{example}

\subsection{Application to Linear Dynamic Equations}\label{secatlde}

In this section, we prove a connection between the Volterra integral equations and initial value problems of the form
\begin{equation}
\begin{cases}
y^{\Delta^{n}}(t)+\displaystyle\sum_{i=0}^{n-1}p_{n-i}(t)y^{\Delta^{i}}(t)=q(t)\quad\text{for}\ t\in\T^{\kappa^{n}}\\
y^{\Delta^{i}}(s)=y_{i}\quad\text{for}\ i=0,1,\ldots,n-1,
\end{cases}\label{atldeeq1}
\end{equation}
where $n\in\N$, $s\in\T^{\kappa^{n-1}}$ and $p_{i},q\in\crd{}(\T,\R)$ for $i=1,2,\ldots,n$.

\begin{lemma}\label{atldelm1}
If a function $y\in\crd{n}(\T,\R)$ solves \eqref{atldeeq1}, then $y^{\Delta^{n}}$ solves
\begin{equation}
\varphi(t)=-\int_{s}^{t}\sum_{i=0}^{n-1}p_{n-i}(t)\func{h}_{n-i-1}\big(t,\sigma(\eta)\big)\varphi(\eta)\Delta\eta+q(t)-\sum_{i=0}^{n-1}\sum_{k=0}^{n-i-1}y_{k+i}p_{n-i}(t)\func{h}_{k}(t,s)\quad\text{for}\ t\in\T^{\kappa^{n}}.\label{atldelm1eq1}
\end{equation}
\end{lemma}

\begin{proof}
Let the function $y\in\crd{n}(\T,\R)$ solve \eqref{atldeeq1}.
From Taylor's formula (\cite[Theorem~1.113]{MR1843232}) we get
\begin{equation}
y^{\Delta^{i}}(t)=\sum_{k=0}^{n-i-1}y_{k+i}\func{h}_{k}(t,s)+\int_{s}^{t}\func{h}_{n-i-1}\big(t,\sigma(\eta)\big)y^{\Delta^{n}}(\eta)\Delta\eta\label{atldelm1prfeq1}
\end{equation}
for all $t\in\T^{\kappa^{i}}$ and $i=0,1,\ldots,n-1$.
Substituting \eqref{atldelm1prfeq1} into \eqref{atldeeq1}, we have
\begin{equation}
y^{\Delta^{n}}(t)+\sum_{i=0}^{n-1}p_{n-i}(t)\bigg(\sum_{k=0}^{n-i-1}y_{k+i}\func{h}_{k}(t,s)+\int_{s}^{t}\func{h}_{n-i-1}\big(t,\sigma(\eta)\big)y^{\Delta^{n}}(\eta)\Delta\eta\bigg)=q(t)\quad\text{for all}\ t\in\T^{\kappa^{n}}\notag
\end{equation}
proving that $y^{\Delta^{n}}$ solves \eqref{atldelm1eq1}.
\end{proof}

It is not hard to see that \eqref{atldelm1eq1} is of the form \eqref{gvieskeq1} with $\lambda=-1$, $\kernel(t,s)=\sum_{i=0}^{n-1}p_{n-i}(t)\func{h}_{n-i-1}\big(t,\sigma(s)\big)$ and $\forcing(t)=q(t)-\sum_{i=0}^{n-1}\sum_{k=0}^{n-i-1}y_{k+i}p_{n-i}(t)\func{h}_{k}(t,s)$ for $t\in\T^{\kappa^{n}}$ and $s\in\T^{\kappa^{n-1}}$.
The uniqueness of the solution $\varphi$ of the initial value problem \eqref{atldeeq1} corresponds exactly to that of \eqref{gvieskeq1}.

\begin{remark}\label{atldermk1}
Method introduced above can be used to solve integro-dynamic equations of the form
\begin{equation}
\begin{cases}
y^{\Delta^{n}}(t)+\displaystyle\sum_{i=0}^{n-1}p_{n-i}(t)y^{\Delta^{i}}(t)+\sum_{j=0}^{m}\int_{a}^{t}\kernel_{j}(t,\eta)y^{\Delta^{j}}(\eta)\Delta\eta=q(t)\quad\text{for}\ t\in\T^{\kappa^{n}}\\
y^{\Delta^{i}}(a)=y_{i}\quad\text{for}\ i=0,1,\ldots,n-1,
\end{cases}\notag
\end{equation}
where $m,n\in\N$, $a\in\T^{\kappa^{n-1}}$ and $p_{i},q\in\crd{}(\T,\R)$ for $i=1,2,\ldots,n$ and $\kernel_{j}\in\crd{}(\T\times\T,\R)$ for $j=0,1,\ldots,m$.
\end{remark}

\subsection{Kernels of Polynomial Type}\label{seckopt}

In this section, we confine our attention to Volterra integral equations with a kernel of the following type
\begin{equation}
\kernel(t,s)=\sum_{i=0}^{n-1}p_{n-i}(t)\func{h}_{n-i-1}\big(t,\sigma(s)\big)\quad\text{for}\ t\in[a,b]_{\T},\label{kopteq1}
\end{equation}
where $n\in\N$ and $\{p_{k}\}_{k=1}^{n}$ are complex-valued rd-continuous functions on $[a,b]_{\T}$.

Let $\lambda\in\R$, and consider the initial value problem
\begin{equation}
\begin{cases}
y^{\Delta^{n}}(t)-\lambda\displaystyle\sum_{i=0}^{n-1}p_{n-i}(t)y^{\Delta^{i}}(t)=0\quad\text{for}\ t\in\T^{\kappa^{n}}\\
y^{\Delta^{i}\sigma}(s)=\delta_{i,n-1}\quad\text{for}\ i=0,1,\ldots,n-1,
\end{cases}\label{kopteq2}
\end{equation}
where $s\in\T^{\kappa^{n-1}}$, $p_{i}\in\crd{}(\T,\R)$ for $i=1,2,\ldots,n$ and $\delta$ is the Kronecker's delta.

The following lemma shows that the resolvent kernel $\Gamma$ of \eqref{gvieskeq1} with a kernel $\kernel$ of the form \eqref{kopteq1} can be computed by means of the unique solution of the dynamic equation \eqref{kopteq2}.

\begin{lemma}\label{koptlm1}
If $\Gamma$ is the resolvent kernel of the kernel $\kernel$ given in \eqref{kopteq1}, then for $\lambda\neq0$, we have
\begin{equation}
\Gamma(\lambda;t,s)=\frac{1}{\lambda}y^{\Delta^{n}}(\lambda;t,s)\quad\text{for}\ t\in\T^{\kappa^{n}}\ \text{and}\ s\in\T^{\kappa^{n-1}},\notag
\end{equation}
where $y(\lambda;\cdot,s)$ denotes the unique solution of \eqref{kopteq2}.
\end{lemma}

\begin{proof}
Then using Lemma~\ref{atldelm1}, we see that the corresponding integral equation for \eqref{kopteq2} is
\begin{equation}
\varphi(t)=\lambda\int_{\sigma(s)}^{t}\sum_{i=0}^{n-1}p_{n-i}(t)\func{h}_{n-i-1}\big(t,\sigma(\eta)\big)\varphi(\eta)\Delta\eta+\lambda\sum_{i=0}^{n-1}p_{n-i}(t)\func{h}_{n-i-1}\big(t,\sigma(s)\big)\notag
\end{equation}
or simply
\begin{equation}
\varphi(t)=\lambda\int_{\sigma(s)}^{t}\kernel(t,\eta)\varphi(\eta)\Delta\eta+\lambda\kernel(t,s)\notag
\end{equation}
by letting $\varphi:=y^{\Delta^{n}}(\lambda;\cdot,s)$ on $\T^{\kappa^{n}}$.
By Definition~\ref{rkdf2}, Theorem~\ref{rkthm1} and Theorem~\ref{rkthm3}, we have 
\begin{align}
\varphi(t)=&\lambda\int_{\sigma(s)}^{t}\Gamma(\lambda;t,\eta)\lambda\kernel(\eta,s)\Delta\eta+\lambda\kernel(t,s)\notag\\
=&\lambda\int_{\sigma(s)}^{t}\bigg(\sum_{\ell=0}^{\infty}\kernel_{\ell}(t,\eta)\lambda^{\ell}\bigg)\lambda\kernel(\eta,s)\Delta\eta+\lambda\kernel(t,s)\notag\\
=&\lambda\sum_{\ell=0}^{\infty}\bigg(\int_{\sigma(s)}^{t}\kernel_{\ell}(t,\eta)\kernel(\eta,s)\Delta\eta\bigg)\lambda^{\ell+1}+\lambda\kernel(t,s)\notag\\
=&\lambda\sum_{\ell=0}^{\infty}\kernel_{\ell+1}(t,s)\lambda^{\ell+1}+\lambda\kernel(t,s)\notag\\
=&\lambda\Gamma(\lambda;t,s)\notag
\end{align}
for all $t\in\T^{\kappa^{n}}$.
Therefore, if \eqref{gvieskeq1} has a kernel of type \eqref{kopteq1}, then the resolvent kernel $\Gamma$ is given by
\begin{equation}
\frac{1}{\lambda}y^{\Delta^{n}}(\lambda;t,s)\quad\text{for}\ t\in\T^{\kappa^{n}}\ \text{and}\ s\in\T^{\kappa^{n-1}},\notag
\end{equation}
where $y(\lambda;\cdot,s)$ is the unique solution of \eqref{kopteq2}.
\end{proof}

As an application of Lemma~\ref{koptlm1}, we have the following illustrative example.

\begin{example}\label{koptex1}
Let us obtain the resolvent kernel for the integral equation
\begin{equation}
\varphi(t)=\int_{a}^{t}\func{h}_{1}\big(t,\sigma(\eta)\big)\varphi(\eta)\Delta\eta+\forcing(t)\quad\text{for}\ t\in[a,b]_{\T},\label{koptex1eq1}
\end{equation}
where $\forcing$ is a rd-continuous function defined on $[a,b]_{\T}$.
Here, we have $\lambda=1$, $\kernel=\func{h}_{1}^{\sigma_{2}}$ on $\Omega(b,a)$, and $p_{1}\equiv0$ and $p_{2}\equiv1$ on $[a,b]_{\T}$.
Then, \eqref{kopteq2} reduces to
\begin{equation}
\begin{cases}
y^{\Delta^{2}}(t)-y(t)=0\quad\text{for}\ t\in\T^{\kappa^{2}}\\
y^{\sigma}(s)=0\ \text{and}\ y^{\Delta\sigma}(s)=1,
\end{cases}\notag
\end{equation}
whose unique solution is given by
\begin{equation}
y(t,s;1)=\frac{1}{2\big(1+\mu(s)\big)}\func{e}_{1}(t,s)+\frac{1}{2\big(1-\mu(s)\big)}\func{e}_{-1}(t,s)\quad\text{for}\ s,t\in\T\notag
\end{equation}
provided that $-1\in\creg{}(\T,\R)$.
The resolvent kernel of \eqref{koptex1eq1} is therefore given by
\begin{equation}
\Gamma(t,s;1)=y^{\Delta^{2}}(t,s;1)=y(t,s;1)=\frac{1}{2\big(1+\mu(s)\big)}\func{e}_{1}(t,s)+\frac{1}{2\big(1-\mu(s)\big)}\func{e}_{-1}(t,s)\notag
\end{equation}
for $t\in\T^{\kappa^{2}}$ and $s\in\T^{\kappa}$.
\end{example}

\subsection{Kernels of Convolution Type}\label{seckoct}

In this section, we shall study Volterra integral equation of the second kind on an unbounded domain.
We start this section with the definition of the generalized Laplace transform on time scales (see \cite{MR1843232,MR2320804,bo/gu/ka10}).

\begin{definition}[Laplace transform]\label{koctdf0}
Let $\sup\T=\infty$.
For a given $f:[a,\infty)_{\T}\to\R$, the Laplace transform is defined by
\begin{equation}
\lt\{f\}(z;a):=\int_{a}^{\infty}\forcing(\eta)\func{e}_{\ominus_{\mu}z}\big(\sigma(\eta),a\big)\Delta\eta\quad\text{for}\ z\in\mathcal{D}_{f},\notag
\end{equation}
where $D_{f}$ consists of complex regressive constants for which the improper integral coverges.
\end{definition}

The following definition is firstly introduced in \cite[Definition~2.1]{MR2320804}.

\begin{definition}[Shift of a function]\label{koctdf1}
For a given $f:[a,\infty)_{\T}\to\R$, the solution of the partial dynamic equation
\begin{equation}
\begin{cases}
\varphi^{\Delta_{1}}(t,s)+\varphi^{\Delta_{2}\sigma_{2}}(t,s)=0,&(t,s)\in[s,\infty)_{\T^{\kappa}}\times[a,\infty)_{\T^{\kappa}}\\
\varphi(t,a)=f(t),&t\in[a,\infty)_{\T}
\end{cases}
\end{equation}
is called the shift of $f$ and is denoted by $\widehat{f}$.
\end{definition}

\begin{definition}[Convolution]\label{koctdf2}
For given functions $f,g:[a,\infty)_{\T}\to\R$, their convolution $f\ast g$ is defined by
\begin{equation}
(f\ast g)(t):=\int_{a}^{t}\widehat{f}\big(t,\sigma(\eta)\big)g(\eta)\Delta\eta\quad\text{for}\ t\in[a,\infty)_{\T}\notag
\end{equation}
or equivalently
\begin{equation}
(f\ast g)(t):=\int_{a}^{t}f(\eta)\widehat{g}\big(t,\sigma(\eta)\big)\Delta\eta\quad\text{for}\ t\in[a,\infty)_{\T}.\notag
\end{equation}
\end{definition}

Let $\sup\T=\infty$, and we consider the following type of Volterra integral equation of the second kind
\begin{equation}
\varphi(t)=\lambda\int_{a}^{t}\widehat{\kernel}\big(t,\sigma(\eta)\big)\varphi(\eta)\Delta\eta+\forcing(t)\quad\text{for}\ t\in[a,\infty)_{\T},\label{kocteq1}
\end{equation}
where $\kernel,\forcing\in\crd{}([a,\infty)_{\T},\R)$ and $\lambda\in\R$.

If $\kernel,\forcing$ are of some exponential order, i.e., there exist $\constl,\constm\in\R_{0}^{+}$ and $\alpha,\beta\in\creg{+}([a,\infty)_{\T},\R)$ such that $|\kernel(t)|\leq\constm\func{e}_{\beta}(t,a)$ and $|\forcing(t)|\leq\constl\func{e}_{\alpha}(t,a)$ for all $t\in[a,\infty)_{\T}$.
Then the Laplace transform of \eqref{kocteq1} exists on $(\gamma,\infty)_{\R}$, where $\gamma:=\max\{\alpha,\beta\}$ (see \cite{bo/gu/ka10}).



We would like to illustrate the discussion above with the following example.

\begin{example}\label{koctex1}
Consider the equation
\begin{equation}
\varphi(t)=2\int_{a}^{t}\cos_{1}\big(t,\sigma(\eta)\big)\varphi(\eta)\Delta\eta+\sin_{1}(t,a)\quad\text{for}\ t\in[a,\infty)_{\T}.\label{koctex1eq1}
\end{equation}
Clearly, when \eqref{koctex1eq1} is compared to \eqref{kocteq1}, we have $\kernel(t)=2\cos_{1}(t,a)$ and $\forcing(t)=\sin_{1}(t,a)$ for $t\in[a,\infty)_{\T}$ .
Since, we have
\begin{equation}
\lt\{\cos_{1}(\cdot,a)\}(z;a)=\frac{z}{z^{2}+1}\quad\text{and}\quad\lt\{\sin_{1}(\cdot,a)\}(z;a)=\frac{1}{z^{2}+1}\quad\text{for}\ z\in(1,\infty)_{\R},\notag
\end{equation}
we get
\begin{equation}
\lt\{\varphi\}(z;a)=\frac{2z}{z^{2}+1}\lt\{\varphi\}(z;a)+\frac{1}{z^{2}+1},\notag
\end{equation}
which yields
\begin{equation}
\lt\{\varphi\}(z;a)=\frac{1}{(z-1)^{2}}.\label{koctex1eq2}
\end{equation}
Taking the inverse Laplace transform of \eqref{koctex1eq2}, we obtain the solution as
\begin{equation}
\varphi(t)=\func{m}_{1}(t,a)\func{e}_{1}(t,a)\quad\text{for}\ t\in[a,b]_{\T},\notag
\end{equation}
where
\begin{equation}
\func{m}_{\lambda}(t,s):=\int_{s}^{t}\frac{1}{1+\lambda\mu(\eta)}\Delta\eta\quad\text{for}\ s,t\in\T\ \text{and}\ \lambda\in\creg{}([a,b]_{\T},\R).\notag
\end{equation}
However, one should (can) justify validity of these solutions since uniqueness of the Laplace transform on time scales
has not been given yet, that is, the inverse may not be unique. 
\end{example}

\begin{remark}\label{koctrmk1}
Let $m,n\in\N$, then the solution of the following system of equations with convolution type kernels
\begin{equation}
\begin{cases}
\varphi_{1}(t)=\displaystyle\sum_{j=1}^{n}\int_{a}^{t}\widehat{\kernel}_{1,j}\big(t,\sigma(\eta)\big)\varphi_{j}(\eta)\Delta\eta+f_{1}(t)\\
\varphi_{2}(t)=\displaystyle\sum_{j=1}^{n}\int_{a}^{t}\widehat{\kernel}_{2,j}\big(t,\sigma(\eta)\big)\varphi_{j}(\eta)\Delta\eta+f_{2}(t)\\
\quad\vdots\\
\varphi_{m}(t)=\displaystyle\sum_{j=1}^{n}\int_{a}^{t}\widehat{\kernel}_{m,j}\big(t,\sigma(\eta)\big)\varphi_{j}(\eta)\Delta\eta+f_{m}(t)
\end{cases}\quad\text{for}\ t\in[a,\infty)_{\T},\notag
\end{equation}
where $\kernel_{i,j},f_{i}\in\crd{}([a,\infty)_{\T},\R)$ for $i=1,2,\ldots,m$ and $j=1,2,\ldots,n$,
can be obtained by making use of the Laplace transform.
\end{remark}

We have another illustrative example below.

\begin{example}\label{koctex2}
Let us solve the following system of integral equations
\begin{equation}
\begin{cases}
\varphi_{1}(t)=-2\displaystyle\int_{a}^{t}\func{e}_{2}\big(t,\sigma(\eta)\big)\varphi_{1}(\eta)\Delta\eta+\int_{a}^{t}\varphi_{2}(\eta)\Delta\eta+1\\
\varphi_{2}(t)=-\displaystyle\int_{a}^{t}\varphi_{1}(\eta)\Delta\eta+4\int_{a}^{t}\func{h}_{1}\big(t,\sigma(\eta)\big)\varphi_{2}(\eta)\Delta\eta+4\func{h}_{1}(t,a)
\end{cases}\quad\text{for}\ t\in[a,\infty)_{\T}.\label{koctex2eq1}
\end{equation}
Taking the Laplace transform of \eqref{koctex2eq1} and denoting by $\Phi_{1}:=\lt\{\varphi_{1}\}(\cdot;a)$ and $\Phi_{2}:=\lt\{\varphi_{2}\}(\cdot;a)$, we get
\begin{equation}
\begin{cases}
\Phi_{1}(z)=\dfrac{2}{z-2}\Phi_{1}(z)+\dfrac{1}{z}\Phi_{2}(z)+\dfrac{1}{z}\\
\Phi_{2}(z)=-\dfrac{1}{z}\Phi_{1}(z)+\dfrac{4}{z^{2}}\Phi_{2}(z)+\dfrac{4}{z^{2}}.
\end{cases}\label{koctex2eq2}
\end{equation}
Solving \eqref{koctex2eq2} for $\Phi_{1}$ and $\Phi_{2}$, we obtain
\begin{equation}
\begin{aligned}
\Phi_{1}(z)=&\dfrac{z}{(z+1)^{2}}=\frac{1}{z+1}-\frac{1}{(z+1)^{2}}\\
\Phi_{2}(z)=&\dfrac{3z+2}{(z-2)(z+1)^{2}}=\frac{8}{9}\frac{1}{z-2}-\frac{8}{9}\frac{1}{z+1}+\frac{1}{3}\frac{1}{(z+1)^{2}}.
\end{aligned}\label{koctex2eq3}
\end{equation}
Taking now the inverse Laplace transform of \eqref{koctex2eq3} to obtain the solution to \eqref{koctex2eq1}, we get
\begin{equation}
\begin{aligned}
\varphi_{1}(t)=&\func{e}_{\ominus1}(t,a)-\func{m}_{\ominus1}(t,a)\func{e}_{\ominus1}(t,a)\\
\varphi_{2}(t)=&\frac{8}{9}\func{e}_{2}(t,a)-\frac{8}{9}\func{e}_{\ominus1}(t,a)+\frac{1}{3}\func{m}_{\ominus1}(t,a)\func{e}_{\ominus1}(t,a)
\end{aligned}\notag
\end{equation}
for $t\in[a,b]_{\T}$.
The validity of these solutions should (can) be justified because of the absence of uniqueness of the Laplace transform on time scales. 
\end{example}

\section{Generalized Volterra Integral Equation of the First Kind}\label{secgviefk}

In this section, we shall consider
\begin{equation}
\int_{a}^{t}\kernel(t,\eta)\varphi(\eta)\Delta\eta=\forcing(t)\quad\text{for}\ t\in[a,b]_{\T},\label{gviefkeq1}
\end{equation}
where $\kernel\in\crd{}(\Omega(a,b),\R)$ and $\forcing\in\crd{}([a,b]_{\T},\R)$.

\subsection{Transform to the Second Kind}\label{sectsk}

In this section, we will write \eqref{gviefkeq1} in the form of \eqref{gvieskeq1} under some conditions.
We see obviously that if $\forcing(a)\neq0$, then \eqref{gviefkeq1} has no solutions.

\begin{theorem}\label{tskthm1}
Assume that $\kernel,\kernel^{\Delta_{1}}\in\crd{}(\Omega(a,b),\R)$, $\forcing\in\crd{1}([a,b],\R)$ and that $\kernel\big(\sigma(t),t\big)\neq0$ for all $t\in[a,b]_{\T^{\kappa}}$.
Then \eqref{gviefkeq1} admits a unique solution on $[a,b]_{\T^{\kappa}}$.
\end{theorem}

\begin{proof}
Differentiating \eqref{gviefkeq1}, we get
\begin{equation}
\int_{a}^{t}\kernel^{\Delta_{1}}(t,\eta)\varphi(\eta)\Delta\eta+\kernel\big(\sigma(t),t\big)\varphi(t)=f^{\Delta}(t)\quad\text{for}\ t\in[a,b]_{\T^{\kappa}}\notag
\end{equation}
or equivalently,
\begin{equation}
\varphi(t)=-\int_{a}^{t}\frac{\kernel^{\Delta_{1}}(t,\eta)}{\kernel\big(\sigma(t),t\big)}\varphi(\eta)\Delta\eta+\frac{f^{\Delta}(t)}{\kernel\big(\sigma(t),t\big)}\quad\text{for}\ t\in[a,b]_{\T^{\kappa}},\notag
\end{equation}
which is a Volterra type integral equation of the first kind.
An application of Theorem~\ref{eauosthm1} completes the proof.
\end{proof}

\begin{remark}\label{tskrmk1}
In the case $\kernel\big(\sigma(t),t\big)=0$ for some $t\in[a,b]_{\T^{\kappa}}$ and $\kernel\in\crd{}(\Omega(a,b),\R)$ has rd-continuous partial derivatives of higher order with respect to its first component. 
Then we may proceed similarly for showing existence and uniqueness of solutions to \eqref{gviefkeq1}.
\end{remark}

We finalize the paper with the following simple application of Theorem~\ref{tskthm1}.

\begin{example}\label{tskex1}
Consider the following equation
\begin{equation}
\int_{a}^{t}\cos_{1}\big(t,\sigma(\eta)\big)\varphi(\eta)\Delta\eta=\func{h}_{1}(t,a)\quad\text{for}\ t\in[a,b]_{\T},\label{tskex1eq1}
\end{equation}
where we have $\kernel=\cos_{1}^{\sigma_{2}}$ on $\Omega(a,b)$, and $f=\func{h}_{1}(\cdot,a)$ on $[a,b]_{\T}$.
Differentiating \eqref{tskex1eq1}, we get
\begin{equation}
\cos_{1}\big(\sigma(t),\sigma(t)\big)\varphi(t)-\int_{a}^{t}\sin_{1}\big(t,\sigma(\eta)\big)\varphi(\eta)\Delta\eta=1\quad\text{for}\ t\in[a,b]_{\T^{\kappa}}\notag
\end{equation}
or equivalently
\begin{equation}
\varphi(t)=\int_{a}^{t}\sin_{1}\big(t,\sigma(\eta)\big)\varphi(\eta)\Delta\eta+1\quad\text{for}\ t\in[a,b]_{\T^{\kappa}}.\label{tskex1eq2}
\end{equation}
Clearly, \eqref{tskex1eq2} involves a kernel of convolution type.
Hence, taking the Laplace transform of \eqref{tskex1eq2}, we have
\begin{equation}
\Phi(z)=\frac{1}{z^{2}+1}\Phi(z)+\frac{1}{z},\notag
\end{equation}
which yields
\begin{equation}
\Phi(z)=\frac{1+z^{2}}{z^{3}}=\frac{1}{z^{3}}+\frac{1}{z}.\label{tskex1eq3}
\end{equation}
Since the Laplace transform of the function $\func{h}_{2}(\cdot,a)+1$ gives the right-hand side of \eqref{tskex1eq3}, we get
\begin{equation}
\varphi(t)=\func{h}_{2}(t,a)+1\quad\text{for}\ t\in[a,b]_{\T},\notag
\end{equation}
which can be shown to be the desired solution to \eqref{tskex1eq2} and also \eqref{tskex1eq1}. 
\end{example}

\section{Further Comments}\label{fc}

Consider the following nonlinear type of Volterra integral equation:
\begin{equation}
\varphi(t)=\lambda\int_{a}^{t}F\big(t,\eta,\varphi(\eta)\big)\Delta\eta+\forcing(t)\quad\text{for}\ t\in[a,b]_{\T}\ \text{and}\ \lambda\in\R,\label{fceq1}
\end{equation}
where $\forcing$ satisfies \ref{h2} and
\begin{enumerate}[label={(H\arabic*)},leftmargin=*,ref=(H\arabic*)]
\setcounter{enumi}{2}
\item\label{h3} $F:\Omega(a,b)\times[-\alpha,\alpha]_{\R}\to\R$ is rd-continuous on $\Omega(a,b)$ and continuous on $[-\alpha,\alpha]_{\R}$ (for some $\alpha\in\R^{+}$).
    There exist $\constl,\constm\in\R^{+}$ such that $F$ is bounded by $\constm$, i.e.,
    \begin{equation}
    |F(t,s,x)|\leq\constm\quad\text{for all}\ (t,s,x)\in\Omega(a,b)\times[-\alpha,\alpha]_{\R}\notag
    \end{equation}
    and satisfies the Lipschitz condition with the constant $\constl$, i.e.,
    \begin{equation}
    |F(t,s,x)-F(t,s,y)|\leq\constl|x-y|\quad\text{for all}\ x,y\in[-\alpha,\alpha]_{\R},\notag
    \end{equation}
    where $(t,s)\in\Omega(a,b)$.
\end{enumerate}
Denoting by $c:=\max[a,a+\delta]_{\T}$, where $\delta:=\min\{b-a,\alpha/\constm\}$,
one can define the sequence of successive approximations $\{\varphi_{n}\}_{n\in\N_{0}}$ by
\begin{equation}
\varphi_{n}(t)=
\begin{cases}
\forcing(t),&n=0\\
\displaystyle\int_{a}^{t}F\big(t,\eta,\varphi_{n-1}(\eta)\big)\Delta\eta,&n\in\N
\end{cases}\quad\text{for}\ t\in[a,c]_{\T}\notag
\end{equation}
as in \eqref{eauosthmprfeq2} to prove that $\{\varphi_{n}\}_{n\in\N_{0}}$ is uniformly convergent on $\Omega(a,c)\times[-\alpha,\alpha]_{\R}$ and
the limiting function is the unique solution of \eqref{fceq1} on $[a,c]_{\T}$, which is extendable to $[a,\sigma(c)]_{\T}$ (see \cite{ka10}).

In the nonlinear case, we can give an error bound for the Picard iterates.
Let $\varphi$ be the unique solution of \eqref{fceq1} and $\{\varphi_{n}\}_{n\in\N_{0}}$ be the Picard iterates defined by
\begin{equation}
\varphi_{n}(t):=\lambda\int_{a}^{t}F\big(t,\eta,\varphi_{n-1}(\eta)\big)\Delta\eta+\forcing(t)\quad\text{for}\ t\in[a,c]_{\T}\ \text{and}\ n\in\N,\notag
\end{equation}
where $\varphi_{0}\in\crd{}([a,b]_{\T},\R)$ is chosen arbitrarily.
Then, by induction, we have
\begin{equation}
|\varphi(t)-\varphi_{n}(t)|\leq|\lambda|^{n+1}\constm\constl^{n}\func{h}_{n+1}(t,a)+|\lambda|^{n-1}\constn\constl^{n}\func{h}_{n}(t,a)\quad\text{for all}\ t\in[a,c]_{\T}\ \text{and}\ n\in\N_{0},\label{fccrl1eq1}
\end{equation}
where
\begin{equation}
\constn:=\sup_{t\in[a,c]_{\T}}|\forcing(t)-\varphi_{0}(t)|.\notag
\end{equation}

In \cite[Theorem~4.1]{MR2193645}, an upper bound for the generalized polynomials is given by
\begin{equation}
\func{h}_{k}(t,s)\leq\frac{(t-s)^{k}}{k!}\quad\text{for all}\ s,t\in\T\ \text{with}\ t\geq s\ \text{and}\ k\in\N_{0}.\notag
\end{equation}
This fact shows that the right-hand side of \eqref{fccrl1eq1} tends to $0$ as $n\to\infty$.

The results of this paper can be generalized to the following system of Volterra integral equations:
\begin{equation}
\pmb{\varphi}(t)=\lambda\int_{a}^{t}\pmb{\kernel}(t,\eta)\pmb{\varphi}(\eta)\Delta\eta+\pmb{\forcing}(t)\quad\text{for}\ t\in[a,b]_{\T}\ \text{and}\ \lambda\in\C,\notag
\end{equation}
where $m\in\N$, $\pmb{\kernel}:\Omega(a,b)\to\R^{m\times m}$ and $\pmb{\forcing}:[a,b]_{\T}\to\R^{m}$ are locally bounded $\Delta$-integrable matrix functions.

Next, as we have mentioned in \S~\ref{secpots}, we now restate and give the complete proof of \cite[Lemma~1]{MR2506155} concerning double-iterated integrals.

\begin{theorem}[Change of Integration Order]\label{ssecithm11}
Let $a\in\T$, $b\in\T$ with $b>a$ and assume that $f:\T\times\T\to\R$ is $\Delta$-integrable on $\{(t,s)\in\T\times\T:\ b>t>s\geq a\}$.
Then,
\begin{equation}
\int_{a}^{b}\int_{a}^{\eta}f(\eta,\xi)\Delta\xi\Delta\eta=\int_{a}^{b}\int_{\sigma(\xi)}^{b}f(\eta,\xi)\Delta\eta\Delta\xi.\notag
\end{equation}
\end{theorem}

\begin{proof}
Define $g:\T\times\T\to\R$ by
\begin{equation}
g(t,s):=
\begin{cases}
f(t,s),&b>t>s\geq a\\
0,&b>t=s\geq a\\
-f(s,t),&b>s>t\geq a
\end{cases}\quad\text{for}\ s,t\in\T.\label{ssecithm11it1eq1}
\end{equation}
Then
\begin{align}
\int_{a}^{b}\int_{a}^{b}g(\eta,\xi)\Delta\xi\Delta\eta=&\int_{a}^{b}\int_{a}^{\eta}g(\eta,\xi)\Delta\xi\Delta\eta+\int_{a}^{b}\int_{\eta}^{b}g(\eta,\xi)\Delta\xi\Delta\eta\notag\\
=&\int_{a}^{b}\int_{a}^{\eta}g(\eta,\xi)\Delta\xi\Delta\eta+\int_{a}^{b}\int_{\sigma(\eta)}^{b}g(\eta,\xi)\Delta\xi\Delta\eta+\int_{a}^{b}\mu(\eta)g(\eta,\eta)\Delta\eta.\notag
\end{align}
Note that the integral variables $\xi$ and $\eta$ in the first and the second integral above vary in $(\eta,\xi)\in\{(t,s)\in\T\times\T:\ b>t>s\geq a\}$ and $(\eta,\xi)\in\{(t,s)\in\T\times\T:\ b>t\geq\sigma(s)\geq a\}\supset\{(t,s)\in\T\times\T:\ b>t>s\geq a\}$, respectively.
And integral of $g$ over the set $\{(t,s)\in\T\times\T:\ b>t\geq\sigma(s)\geq a\}\backslash\{(t,s)\in\T\times\T:\ b>t>s\geq a\}\subset\{(t,t)\in\T\times\T:\ b>t>a\}$ is $0$.
This gives us
\begin{equation}
\int_{a}^{b}\int_{a}^{b}g(\eta,\xi)\Delta\xi\Delta\eta=\int_{a}^{b}\int_{a}^{\eta}f(\eta,\xi)\Delta\xi\Delta\eta-\int_{a}^{b}\int_{\sigma(\eta)}^{b}f(\xi,\eta)\Delta\xi\Delta\eta.\label{ssecithm11it1eq2}
\end{equation}
On the other hand, the function $g$ defined by \eqref{ssecithm11it1eq1} is anti-symmetric on the domain $[a,b)_{\T}\times[a,b)_{\T}$,
and thus the right-hand side of \eqref{ssecithm11it1eq2} is $0$,
which together completes the proof.
\end{proof}

\begin{remark}
The conclusion of Theorem~\ref{ssecithm11} is true when $f\in\crd{}(\Omega(a,b),\R)$.
\end{remark}

\begin{remark}
An $n$-fold integral can be considered as a finite repetition of double-iterated integrals
and Theorem~\ref{ssecithm11} can be applied from the inmost integration to the outer one as many as needed.
\end{remark}

\end{document}